

\magnification\magstep1 
\baselineskip14pt 
\vsize23.5truecm 
\font\bigbf=cmbx12 
\font\csc=cmcsc10 
\font\smallrm=cmr10 at 10truept 
\font\cyr=wncyr10 
\def\j3{{\rm\u{\cyr i}}}

\def\hop{\smallskip\noindent} 
 
\def\Beta{{\rm Beta}}
\def\Dir{{\rm Dir}}
\def\E{{\rm E}} 
 
\def\Var{{\rm Var}} 
\def\d{{\rm d}} 
\def\data{{\rm data}}
\def\info{{\rm info}} 
\def\midd{\,|\,}  
\def\arr{\rightarrow} 
\def\hatt{\widehat} 
\def\tilda{\widetilde}  
\def\dell{\partial} 
\def\sumin{\sum_{i=1}^n}
\def\prodin{\prod_{i=1}^n}

\def\eps{\varepsilon} 
\def\half{\hbox{$1\over2$}}

\def\section{\bigskip} 
\def\subsection{\medskip} 
\def\GDir{{\rm GD}}
\def\true{{\rm true}}
\def\square{{\ \vrule height0.5em width0.5em depth-0.0em}}

\def\fermat#1{\setbox0=\vtop{\hsize4.00pc
        \smallrm\raggedright\noindent\baselineskip9pt
        \rightskip=0.5pc plus 1.5pc #1}\leavevmode
        \vadjust{\dimen0=\dp0
        \kern-\ht0\hbox{\kern-4.00pc\box0}\kern-\dimen0}}

\centerline{\bigbf Bayesian analysis for a generalised 
   Dirichlet process prior}  
\medskip
\centerline{\bf Nils Lid Hjort}
\smallskip
\centerline{\bf University of Oslo} 
\smallskip 
\centerline{\sl November 2000} 

{\bigskip\narrower\noindent
{\csc Abstract.} 
A family of random probabilities is defined and studied.
This family contains the Dirichlet process as a special case, 
corresponding to an inner point in the appropriate parameter space. 
The extension makes it possible to have random means with 
larger or smaller skewnesses as compared to skewnesses  
under the Dirichlet prior, and also in other ways amounts to 
additional modelling flexibility. 
The usefulness of such random probabilities for
use in nonparametric Bayesian statistics is discussed.
The posterior distribution is complicated, but inference
can nevertheless be carried out via simulation, and some exact
formulae are derived for the case of random means. 
The class of nonparametric priors provides an instructive example
where the speed with which the posterior forgets its prior 
with increasing data sample size depends on special aspects of the prior, 
which is a different situation from that of parametric inference. 
\smallskip}

{\smallskip\narrower\noindent 
{\csc Key words:} \sl 
consistency, 
Dirichlet process,
jump sizes, 
nonparametric Bayes, 
random means,
speed of memory loss, 
stochastic equation
\smallskip} 

\section
\centerline{\bf 1. Introduction and summary} 

\hop
The Dirichlet process, introduced in Ferguson (1973, 1974),
continues to be a cornerstone of nonparametric Bayesian statistics, 
where it may be used as a prior for an unknown probability distribution
for data. Various generalisations have been proposed and 
investigated in the literature, making the Dirichlet a special
favourite case of P\'olya trees, of Beta processes, 
of neutral to the right and of tailfree processes 
and of Dirichlet mixtures; see Walker, Damien, Laud and Smith (1998) 
and Hjort (2001) for recent overviews. The purpose of this article 
is to provide yet another generalisation of the Dirichlet 
process and to study some of its properties. 

Write $P\sim\Dir(b,P_0)$ to signify that $P$ is a Dirichlet process 
with parameters $(b,P_0)$ on some sample space $\Omega$, where 
$b$ is positive and $P_0$ a probability distribution;
in particular, for each set $A$, the random probability $P(A)$ 
has a Beta distribution with parameters $(bP_0(A),b\{1-P_0(A)\})$. 
For a function $g$ of interest, consider the random mean 
$\theta=\E\,\{g(X)\midd P\}=\int g\,\d P$. Its moments 
are most conveniently given in terms of a variable $Y=g(\xi)$ 
with $\xi$ being drawn from $P_0$. Then 
$$\E\,\theta=\theta_0, \quad
  \Var\,\theta={1\over 1+b}\E_0\,(Y-\theta_0)^2, \quad 
  \E\,(\theta-\theta_0)^3={2\over (1+b)(2+b)}\E_0\,(Y-\theta_0)^3, \eqno(1.1)$$ 
where $\theta_0=\E_0\,Y$; here `$\E_0$' signifies expectations with
respect to the base distribution $Q_0=P_0g^{-1}$ for $Y$. 
The two first results are in Ferguson (1973) while the third 
may be proved using similar arguments; in Section 3 below 
we actually give a method for finding the full moment sequence. 
For $g$ the indicator function of a set $A$, (1.1) specialises to 
$$\E\,P(A)=p_0, \quad 
  \Var\,P(A)={p_0(1-p_0)\over 1+b}
  \quad {\rm and} \quad 
  \E\,\{P(A)-p_0\}^3={2p_0(1-p_0)(1-2p_0)\over (1+b)(2+b)}, $$
where $p_0=P_0(A)$, agreeing by necessity with moment calculations
from the Beta distribution $(bp_0,b-bp_0)$. 
These equations make clear that there is a good amount of
modelling flexibility with the Dirichlet prior, as one may 
centre the random $P$ at any chosen prior mean distribution $P_0$
and attune $b$ to a desired level of variability. One is then stuck 
with the consequences, however, regarding all further aspects
of the prior, such as the implied skewnesses of random probabilities 
and of random means. 

The generalised class of nonparametric priors to be worked with below
makes it possible to adjust to further aspects of prior knowledge,
for example regarding the skewnesses of random means. 
Let $B_1,B_2,\ldots$ be independent from a suitable distribution on $(0,1)$,
and define random probabilities 
$\gamma_1=B_1$, $\gamma_2=\bar B_1B_2$, $\gamma_3=\bar B_1\bar B_2B_3$
and so on, where $\bar B_i=1-B_i$. Here 
$1-\sum_{j=1}^n\gamma_i=\bar B_1\cdots\bar B_n$, making it easy
to show that the $\gamma_j$s sum to 1 with probability 1. 
This allows us the possibility of defining a random probability 
$$P=\sum_{j=1}^\infty\gamma_j\delta(\xi_j), 
  \quad {\rm with\ independently\ drawn\ }\xi_j{\rm s\ from\ }P_0. \eqno(1.2)$$
Here $\delta(\xi)$ indicates the measure giving full mass 1 to 
location $\xi$. In other words, $P(A)$ for a set $A$ can be represented 
as a random sum of random probabilities $\sum_{\xi_j\in A}\gamma_j$.
As demonstrated in Sethuraman and Tiwari (1982), the Dirichlet
process can also be represented in such a form, corresponding
to the particular case where the distribution for $B_j$ is 
chosen as the $\Beta(1,b)$; see also Sethuraman (1994). 

The prior process given in (1.2) is described by two distributions,
the prior mean $P_0$ and the distribution $H$ on $(0,1)$ 
governing the $B_j$s. Below some properties of this general prior,
say $\GDir(H,P_0)$, are investigated. An attractive class of priors 
emerges by allowing $B_j$ the $\Beta(a,b)$ distribution. We write 
$P\sim\GDir(a,b,P_0)$ to indicate this particular extension of the Dirichlet, 
which corresponds to having $a=1$. It is important to note that 
the Dirichlet becomes an `inner point' in the enlarged class,
in contrast to what is the case for some other proposals,
where the Dirichlet is at `a corner' of the underlying 
parameter space; see Hjort and Ongaro (2001) for some constructions
of that type. 

In Section 2 we demonstrate that the prior (1.2) has large support
in the space of probability measures on the sample space, indicating
that these priors are genuinely nonparametric. Section 3 deals
with Markov chain Monte Carlo simulation methods for $P$ and 
for random means thereof, and uses stochastic identities to 
derive formulae for the full moment sequence of such random means. 
This is used in Section 4 to show that the added flexibility afforded 
by the $\GDir(a,b,P_0)$ class, and a fortiori the more general
$\GDir(H,P_0)$ class, indeed allows skewnesses a larger range 
than that dictated by the Dirichlet. In Section 5 the posterior 
mean estimator for $P$ is exhibited, leading also to explicit
expressions for Bayes estimators of mean parameters. These
are convex combinations of prior means and sample averages. 
This is supplemented in Section 6 with posterior variance 
formulae, which are used in Section 7 to show that not only 
the posterior mean estimator but also the posterior distribution as such 
becomes consistent, in the sense of being able to recapture
any true distribution underlying the data, as the sample size
increases. Interestingly we find that the speed with which this
convergence takes place depends on aspects of the prior process;
in particular, this speed is sometimes faster and sometimes 
slower than the rate $O(n^{-1})$ found for all models with a 
finite number of parameters. The posterior process itself is
somewhat complicated. It is exhibited in Section 8. Then in Section 9
some results for distributions of random means are given 
before we offer a list of concluding remarks in Section 10. 

\section
\centerline{\bf 2. Large supports} 

\hop 
In parametric Bayesian statistics a prior is in effect 
placed on some set of densities, say on 
${\cal M}_0=\{f_\alpha(\cdot)\colon\alpha\in R\}$,
where $R$ is a finite-dimensional parameter set
indexing the $f_\alpha$ densities on the sample space $\Omega$. 
But this set is a thin one in the space $\cal M$ of all
distributions on $\Omega$, topologically speaking; 
natural neighbourhoods around given distributions 
are given zero prior probability. 
The situation is different for genuinely nonparametric priors,
as demonstrated for example by Ferguson (1973) for the 
Dirichlet process prior. We show here that the generalised 
priors (1.2) continue to have such large supports in the 
space $\cal M$. The support is the set of probability 
distributions $Q$ such that every neighbourhood around it
has positive probability. 

In this section we take the distribution $H$ of $B_j$s 
to have full support $[0,1]$. The key observation is that 
if only $P_0(A)>0$, then the possible outcomes of 
$P(A)=\sum_{\xi_j\in A}\gamma_j$ fill out all of $(0,1)$.
More generally, if $A_1,\ldots,A_k$ are disjoint sets with 
positive $P_0$ probability and union probability less than 1,
then the distribution of $(P(A_1),\ldots,P(A_k))$ 
has a density being positive on the $k$-simplex of 
$(p_1,\ldots,p_k)$ with positive components and sum less than~1. 

\subsection
{\csc 2.1. Support under strong convergence.} 
Strong or set-wise convergence $Q_n\arr Q$, for probability
measures on $\Omega$, means that $Q_n(A)\arr Q(A)$ for all sets $A$. 
A basis for neighbourhoods around a given $Q$ is the class of 
$$U=U(Q;A_1,\ldots,A_m,\eps_1,\ldots,\eps_m)
   =\{Q'\colon |Q'(A_j)-Q(A_j)|<\eps_j 
    {\rm\ for\ }j=1,\ldots,m\}, \eqno(2.1)$$
where $m$ is any integer, the $A_j$s are measurable subsets, 
and the $\eps_j$s positive. Here 
$${\rm supp}\{\GDir(H,P_0)\}=\{Q\colon Q<<P_0\}, \eqno(2.2)$$ 
the set of measures absolutely continuous with respect to $P_0$. 
To see this, let $Q<<P_0$. It suffices to show that $U$ has 
positive probability when the $A_j$s form a measurable partition.
If $P_0(A_j)=0$ then $Q(A_j)=0$ and $P(A_j)=0$ too.
Hence $U$ has positive probability if $U'$ has, where 
$U'=\{Q'\colon |Q'(A_j)-Q(A_j)|<\eps_j{\rm\ for\ }j=i_1,\ldots,i_k\}$,
where these are the indexes for which $P_0(A_j)>0$. 
But it follows from the comment made above that 
$\GDir(H,P_0)$ gives positive probability to this event.
Hence $Q$ is in the support. If on the other hand $Q$ is such that
$Q(A)>0$ but $P_0(A)=0$ for some $A$. Then $P(A)=0$ a.s.,
and $\{P\colon |P(A)-Q(A)|<\half Q(A)\}$ does not have positive probability. 

\subsection
{\csc 2.2. Support under weak convergence.} 
Assume now that the sample space has a metric and study 
the topology determined by  weak convergence, where $Q_n\arr Q$ 
means convergence in distribution. A basis for neighbourhoods 
under this topology is the class of (2.1) type sets, 
but with the restriction that the $A_j$ sets are $Q$-continuous, 
that is, $Q(\dell A_j)=0$, where $\dell A_j$ is the boundary set of $A_j$. 
Here 
$${\rm supp}\{\GDir(H,P_0)\}=\{Q\colon
   {\rm supp}(Q)\subset{\rm supp}(P_0)\}. \eqno(2.3)$$ 
Let $Q$ have a support contained in the support of $P_0$. 
In general $P_0(A)=0$ does not imply $Q(A)=0$, but 
this is seen to hold when the set $A$ is $Q$-continuous. 
Hence the arguments used to prove (2.2) can be used with small modifications
to prove (2.3).

\section
\centerline{\bf 3. Stochastic equations, MCMC and random means} 

\hop
In this section a fruitful stochastic equation is exhibited which 
characterises the $\GDir(H,P_0)$ prior process. This is used 
to give a Markov chain Monte Carlo method for simulating 
realisations of the processes and to derive results about random means. 
If one only needs simulated realisations for one or more
random means a simpler Monte Carlo Markov chain suffices. 

\subsection
{\csc 3.1. Stochastic equations and MCMC simulation.} 
Let $P$ have the prior given in (1.2), with a general distribution $H$
for the $B_j$s. Then 
$$\eqalign{
P&=B_1\delta(\xi_1)+\bar B_1B_2\delta(\xi_2)
   +\bar B_1\bar B_2B_3\delta(\xi_3)+\cdots \cr
&=B_1\delta(\xi_1)+\bar B_1\{B_2\delta(\xi_2)+\bar B_2B_3\delta(\xi_3)
   +\bar B_2\bar B_3B_4\delta(\xi_4)+\cdots\} \cr
&=B_1\delta(\xi_1)+\bar B_1P', \cr}$$
where $P'$ is constructed in the very same manner $P$. 
Letting `$=_d$' mean equality in distribution there is accordingly 
a stochastic equation
$$P=_d B\delta(\xi)+\bar BP, \eqno(3.1)$$
where on the right hand side $B$, $\xi$, $P$ are independent,
with $B\sim H$ and $\xi\sim P_0$. One may show that this identity
fully characterises the distribution of $P$. 
Applying (3.1) to a random mean functional $\theta=\int g\,\d P$
one finds that this variable, which may also be expressed as  
$\sum_{j=1}^\infty\gamma_jg(\xi_j)$, satisfies the stochastic equation
$$\theta=_d BY+\bar B\theta, 
  \quad {\rm where\ }Y=g(\xi)\sim P_0g^{-1}
  {\rm\ and\ }B\sim H. \eqno(3.2)$$

A Markov chain $P_1,P_2,\ldots$ may be constructed 
in the space of probability measures on the sample space via 
$$P_n=B_n\delta(\xi_n)+\bar B_nP_{n-1}, $$
where $(B_n,\xi_n)$ are independent copies of $(B,\xi)$.
With arguments parallelling those in Feigin and Tweedie (1989) 
the equilibrium distribution for the chain may be shown to be exactly 
that of our prior process (1.2). 
For a random mean functional, the Markov chain scheme becomes 
$\theta_n=B_nY_n+\bar B_n\theta_{n-1}$, 
with the distribution of $\theta$ under (1.2) as its equilibrium.
This is akin to similar simulation strategies for means of the 
Dirichlet process, worked with in Feigin and Tweedie (1989), 
Guglielmi and Tweedie (2000) and Guglielmi, Holmes and Walker (2001).
See also Paulsen and Hove (1999) for precise results about 
speed of convergence and quality of approximation to the real
distribution with the empirical one observed from simulations. 

Note that when interest lies in one or more random means 
the simpler simulation scheme suffices, as there is no need 
for the full process $P$. We also point out that the moment-correcting 
methods used in Hjort and Ongaro (2000) apply here too and amount to ways 
of easily improving the simulation-based approximations 
of Paulsen and Hove (1999), Guglielmi and Tweedie (2000)
and Guglielmi, Holmes and Walker (2001). The key is that 
the full moment sequence may be uncovered, as we demonstrate next. 

\subsection
{\csc 3.2. Finding the moments.} 
A recursive method of finding all moments for such a $\theta$,
in terms of moments for the null distribution $Q_0$ for $Y=g(\xi)$, 
emerges by writing 
$$(\theta-x)^p=_d \sum_{j=0}^p{p\choose j}B^{p-j}(Y-x)^{p-j}
   \bar B^j(\theta-x)^j, $$
which implies 
$$\E\,(\theta-x)^p={1\over 1-\E\,\bar B^p}
  \sum_{j=0}^{p-1}\E\,B^{p-j}\bar B^j
  \,\E_0\,(Y-x)^{p-j}\,\E\,(\theta-x)^j. \eqno(3.3)$$
This is valid for all $p\ge 1$ for which $\E\,|Y|^p$ is finite,
and for all $x$. One finds in particular $\E\,\theta=\theta_0=\E_0\,Y$ and
$$\eqalign{
\Var\,\theta&={\E\, B^2\over 1-\E\,\bar B^2}\sigma_0^2, \cr
\E\,(\theta-\theta_0)^3
   &={\E\, B^3\over 1-\E\,\bar B^3}\E_0(Y-\theta_0)^3, \cr
\E\,(\theta-\theta_0)^4
   &={1\over 1-\E\,\bar B^4}\Bigl\{\E\, B^4\,\E_0(Y-\theta_0)^4
   +6\,\E\, B^2\bar B^2{\E\, B^2\over 1-\E\,\bar B^2}
      \sigma_0^4 \Bigr\}, \cr} \eqno(3.4)$$
in terms of $\sigma_0^2=\E_0(Y-\theta_0)^2$. 
Further formulae for centralised moments follow from (3.3),
expressed in terms of 
$$M_{i,j}=\E\, B^i\bar B^j=\int_0^1 s^i(1-s)^j\,\d H(s). \eqno(3.5)$$
With $g$ the indicator of a set $A$, the $\theta$ becomes
the random probability $P(A)$, for which we therefore have found 
$\E\, P(A)=P_0(A)=p_0$ and 
$$\Var\,P(A)={M_{2,0}\over 1-M_{0,2}}p_0(1-p_0)
  \quad {\rm and} \quad 
  \E\,\{P(A)-p_0\}^3={M_{3,0}\over 1-M_{0,3}}p_0(1-p_0)(1-2p_0). $$

\section
\centerline{\bf 4. Skewness factors and added flexibility} 

\hop
In this section the increased flexibility of the nonparametric 
prior class is discussed in relation to the skewness of random means. 

\subsection
{\csc 4.1. Skewnesses under the $\GDir(a,b,P_0)$ prior.} 
When $H$ is the $\Beta(a,b)$, 
$$M_{i,j}=\E\, B^i\bar B^j={\Gamma(a+b)\over \Gamma(a)\Gamma(b)}
  {\Gamma(a+i)\Gamma(b+j)\over \Gamma(a+b+i+j)}
  ={a^{[i]}b^{[j]}\over (a+b)^{[i+j]}} $$
for integers $i,j$, where $x^{[i]}=x(x+1)\cdots(x+i-1)=\Gamma(x+i)/\Gamma(x)$. 
This may be used to derive moment formulae under the $\GDir(a,b,P_0)$ prior.
From (3.4) one finds 
$$\Var\,\theta={M_{2,0}\over 1-M_{0,2}}\sigma_0^2=
  {a+1\over a+2b+1}\sigma_0^2={\sigma_0^2\over 1+b^*},
  \quad {\rm with\ }b^*={2b\over 1+a}, \eqno(4.1)$$
and similarly $\Var\,P(A)=P_0(A)\{1-P_0(A)\}/(1+b^*)$. Furthermore, 
$$\eqalign{\E\,(\theta-\theta_0)^3
&={a(a+1)(a+2)/\{(a+b)(a+b+1)(a+b+2)\}
   \over 1 - b(b+1)(b+2)/\{(a+b)(a+b+1)(a+b+2)\}}\E_0(Y-\theta_0)^3 \cr
&={(a+1)(a+2)\over a^2+3a(b+1)+3b^2+6b+2}\E_0(Y-\theta_0)^3. \cr}$$ 
The Dirichlet case is $a=1$ for which the skewness factor is 
$2/\{(b+1)(b+2)\}$; cf.~(1.1). 

Assume a $\Dir(b_0,P_0)$ has been chosen, and consider using 
a more general $\GDir(a,b,P_0)$ instead; using the same base measure 
$P_0$ automatically ensures that the expected values of any random mean
are being matched for the two priors. We may also precisely match 
all pairs of variances of random means through putting $2b/(1+a)=b_0$, 
compare (4.1) and (1.1). This amounts to $a=2b/b_0-1=2x-1$ 
as a function of $x=b/b_0$; notice that $x>\half$, or $b>\half b_0$, 
is required. We may then study the skewness of $\theta$ under 
the Dirichlet versus its value under the $\GDir(a,b,P_0)$. 
The ratio of skewnesses becomes 
$$\eqalign{\rho(x)
&={(a+1)(a+2)\over a^2-1+3a(b+1)+3(b+1)^2}\Big/
    {2\over (b_0+1)(b_0+2)} \cr 
&={2x(2x+1)\over (2x-1)^2-1+3(2x-1)(b_0x+1)+3(b_0x+1)^2}
  {(b_0+1)(b_0+2)\over 2}. \cr}$$
This is a decreasing function, starting for $b=\half b_0$ 
with ratio value $\rho_{\max}$ and ending for large $b$ 
with ratio value $\rho_{\min}$, say, where 
$$\rho_{\max}={(b_0+1)(b_0+2)\over 2+3b_0+(3/4)b_0^2}
  \quad {\rm and} \quad 
  \rho_{\min}={2(b_0+1)(b_0+2)\over 4+6b_0+3b_0^2}. $$
This interval contains the value 1 as an inner point, 
corresponding to $b=b_0$ and $a=1$, the Dirichlet case,
and otherwise portrays the added flexibility through the 
additional $a$ parameter. With $b<b_0$ and accompanying $a=2b/b_0-1$,
the $\GDir(a,b,P_0)$ prior leads to skewnesses bigger in absolute size 
for all random means than with the correspondingly matched Dirichlet
prior; similarly, with $b>b_0$ the skewnesses are reduced in absolute size. 
The interval of skewness ratios stretches from $2/3$ to $4/3$ 
when $b_0$ becomes large. 


\subsection
{\csc 4.2. More flexible kurtosis.}  
A similar exercise is to compute $\E\,(\theta-\theta_0)^4$,
first under the Dirichlet prior $(b_0,P_0)$, and compare 
it with the value obtained under the $\GDir(a,b,P_0)$ process,
having fixed $2b/(1+a)=b_0$ to have the same mean and same
variance. This gives a suitable kurtosis ratio curve $\kappa(x)/\kappa(1)$ 
to study, where $\kappa(x)$ is $\E\,(\theta-\theta_0)^4$ 
computed with $a=2x-1$ and $b=b_0x$, for $x>\half$. 
This ratio curve depends on $b_0$ and the underlying population
kurtosis $\E_0(Y-\theta_0)/\sigma_0^4-3$. 
Carrying out this exercise one finds that the kurtosis 
is larger than under the Dirichlet for $a<1$ and smaller than under 
the Dirichlet for $a>1$. The ratio interval spans for each $b_0$
a reasonable interval containing 1. 

\section
\centerline{\bf 5. Marginal distributions and posterior means} 

\hop
Conditionally on the random $P$, let $X_1,\ldots,X_n$ 
be independently sampled from $P$ in the sample space $\Omega$. 
In this section we consider the marginal distribution
of data and go on to a formula for the predictive distribution,
that is, the posterior mean of $P$. 

\subsection
{\csc 5.1. Marginal distributions.} 
The simultaneous distribution of $P$ and the random sample is given by 
$$\Pr\{P\in C,X_1\in A_1,\ldots,X_n\in A_n\}
  =\E\,\,I\{P\in C\}P(A_1)\cdots P(A_n), \eqno(5.1)$$
required to hold for measurable subsets $C$ of the space of 
probability measures on the space and for all measurable subsets $A_i$; 
see e.g.~Ferguson (1973). Here $I$ denotes an indicator function. 
In particular, 
$$\Pr\{X_1\in A_1,\ldots,X_n\in A_n\}=\E\,P(A_1)\cdots P(A_n). $$
For $n=1$ one finds 
$$\Pr\{X_i\in A\}=\E\, P(A)=P_0(A), $$
adding significance to the interpretation of $P_0$
as the marginal distribution of a single observation.
For $n=2$, 
$$\eqalign{\E\, P(A)P(B)
&=\E\,\sum_{j,k}\gamma_j\gamma_kI\{\xi_j\in A,\xi_k\in B\} \cr 
&=\sum_j\E\,\gamma_j^2\,P_0(A\cap B)
   +2\sum_{j<k}\E\,\gamma_j\gamma_k\,P_0(A)P_0(B) \cr 
&=(1-a_2)P_0(A\cap B)+a_2P_0(A)P_0(B), \cr}$$
where $a_2=\Pr\{X_2\not=X_1\}=1-M_{2,0}/(1-M_{0,2})=2M_{1,1}/(1-M_{0,2})$. 

\eject 
The following identity will be useful. 

{\smallskip\sl
{\csc Lemma.} Let $P$ come from the $\GDir(H,P_0)$ prior (1.2)
and let $A_1,\ldots,A_n$ be disjoint sets. Then 
$$\E\, P(A_1)\cdots P(A_n)=n!\prod_{j=1}^{n-1}{M_{1,j}\over 1-M_{0,j+1}}
  P_0(A_1)\cdots P_0(A_n) \eqno(5.2)$$
for all $n\ge 2$, in terms of the product moments $M_{i,j}$ of (3.5). 
\smallskip} 

{\csc Proof.} 
The (1.2) definition leads to the formula 
$$\E\, P(A_1)\cdots P(A_n)=P_0(A_1)\cdots P_0(A_n)\, 
   n!\sum_{i_1<\cdots<i_n}\E\, \gamma_{i_1}\cdots\gamma_{i_n}, $$
which indeed may be worked further by careful cataloguing of 
$B_k$ and $\bar B_k$ factors entering the product of $\gamma_j$s;
one finds for example that $\E\, P(A_1)P(A_2)=a_2\,P_0(A_1)P_0(A_2)$
with the $a_2$ exhibited above. It is easier, however, 
to use the simultaneous stochastic equations 
$$P(A_i)=_d BY_i+\bar BP(A_i)
  \quad {\rm for\ }i=1,\ldots,n, $$
where $Y_i=I\{\xi\in A_i\}$ for a $\xi$ drawn from $P_0$,
independently of $B\sim H$. That these equations hold simultaneously 
follows from (3.1). All products of two or more $Y_j$s vanish
since the $A_j$s are disjoint. This simplifies the structure of 
$$\prodin P(A_i)=_d\prodin\{BY_i+\bar BP(A_i)\}
   =\sumin BY_i\bar B^{n-1}\theta_{(i)}+\bar B^n\prodin P(A_i), $$
writing $\theta_{(i)}$ for the product of those $n-1$ terms $P(A_j)$ 
for which $j\not=i$. Hence 
$$\E\,\prodin P(A_i)={1\over 1-M_{0,n}}
   \sumin M_{1,n-1}P_0(A_i)\,\E\,\theta_{(i)}. $$
This may now be used to demonstrate (5.2) by induction on $n$,
noting that the formula was seen to hold for $n=2$ above. \square 

\subsection
{\csc 5.2. The posterior mean.} 
For the next development, define 
$$w_n=(n+1){M_{1,n}\over 1-M_{0,n+1}}
     =(n+1){\E\, B\bar B^n\over 1-\E\,\bar B^{n+1}}=(n+1)\delta_n. \eqno(5.3)$$
We take the sample space to be a metric space,
for example a subset of any Euclidean space, where 
we condition on the information in a data point $X_i=x_i$ 
via conditioning on the information $X_i\in S(x_i,\eps)$,
say, an $\eps$-neighbourhood around $x_i$, and then letting $\eps\arr0$. 
For an observed sample, consider the predictive distribution 
$\hatt P$ given by $\hatt P(A)=\E\,\{P(A)\midd\data\}$;
this is also the Bayes estimator of $P$ under squared error loss. 

{\smallskip\sl
{\csc Proposition.} 
Let $P$ follow the prior process $\GDir(H,P_0)$ with an 
atom-free prior mean measure $P_0$, and assume data points 
$X_1=x_1,\ldots,X_n=x_n$ have been observed, and that 
these are distinct. Then the Bayes estimator of $P$ can be represented as 
$$\hatt P=w_nP_0+(1-w_n){1\over n}\sumin\delta(x_i)
    =w_nP_0 + (1-w_n)P_n, \eqno(5.4)$$
a convex combination of prior mean and the empirical distribution.
\smallskip}

{\csc Proof.} 
From (5.1) one may show that 
$$\E\,\{\psi(P)\midd X_1\in A_1,\ldots,X_n\in A_n\}
  ={\E\,\psi(P)P(A_1)\cdots P(A_n)\over \E\, P(A_1)\cdots P(A_n)} $$
for all bounded measurable functions $\psi$, provided the $A_i$s have 
positive $P_0$ measure. In particular, therefore, 
$$\E\,\{P(A)\midd X_1\in A_1,\ldots,X_n\in A_n\}
  ={\E\, P(A)P(A_1)\cdots P(A_n)\over \E\, P(A_1)\cdots P(A_n)}. $$
Take first a set $A$ not meeting the data, which means that 
it is outside the union of data windows $A_i=S(X_i,\eps)$
for small enough $\eps$. 
Then the above quotient, by the use of the lemma, reduces to 
$$\E\,\{P(A)\midd\info_\eps\}=(n+1)\{M_{1,n}/(1-M_{0,n+1})\}\,P_0(A), $$
where $\info_\eps$ signifies the information $X_i\in S(x_i,\eps)$
for $i=1,\ldots,n$. Since the answer is independent of $\eps$,
the probability measure $\hatt P$ must be equal to $w_nP_0$ 
on $\Omega-\{x_1,\ldots,x_n\}$, that is, outside the data values. 

Being a probability measure it must distribute its remaining mass
$1-w_n$ on the $n$ data values $x_1,\ldots,x_n$. 
With these being distinct there must be full symmetry, 
and $\hatt P$ must assign value $(1-w_n)/n$ to each of these.
This proves assertion (5.4). \square 

\subsection
{\csc 5.3. Nonparametric Bayes estimation of means.} 
Consider Bayesian estimation of a random mean $\theta=\int g\,\d P$. 
Under squared error loss 
and with the $\GDir(H,P_0)$ prior, for which the prior guess
is $\theta_0=\int g\,\d P_0$, the estimator is 
$$\hatt\theta=\E\,(\theta\midd\data)=\int g\,\d\hatt P
  =w_n\theta_0+(1-w_n)\bar g_n
  \quad {\rm with\ }\bar g_n=n^{-1}\sumin g(x_i). \eqno(5.5)$$
This follows from (5.4), again under the assumption on there being
no ties in data. With a little more formality, this concerns  
$$\E\,(\theta\midd\data)=\int_{{\cal M}}\theta(P)\,{\cal P}(\d P\midd\data), $$
where ${\cal P}(\cdot\midd\data)$ is the posterior distribution 
on the space $\cal M$ of probability measures on the sample space,
and an ingredient is existence and measurability of $\theta=\theta(P)$. 
A more careful argument, therefore, starts with $g$ equal to a simple function,
a linear combination of indicator functions. For such a $g$
the result follows directly from (5.4). Then pass to the limit
via monotone convergence to make formula (5.5) valid 
for all $g$ for which $\int|g|\,\d P_0$ is finite. 
Measurability comes from it being a limit of linear combinations
of $P(A)$ variables, and existence is guaranteed under the 
minimal condition $\int\log(1+|g|)\,\d P_0$, see Hjort and Ongaro (2000). 

As a special case, when an unknown distribution function $F$ 
of one-dimensional data is to be estimated, the Bayes estimator
takes the form $\hatt F(t)=w_nF_0(t)+(1-w_n)F_n(t)$,
where $F_0$ is the distribution function of $P_0$ and 
$F_n$ is the empirical distribution function. 

\smallskip
{\csc Remark.} 
The Dirichlet case corresponds to a $\Beta(1,b)$ distribution 
for the $B_j$s, and a little algebra on (5.3) shows that for this case 
$w_n=b/(b+n)$. This is a well-known formula for the weight 
a posterior Dirichlet distribution still attaches to its prior, 
also lending strength to the `prior sample size'
interpretation of the $b$ parameter. More nuances are at play
for the general $\GDir(H,P_0)$ case, however, 
as shown in section 7. \square 

\bigskip
\centerline{\bf 6. Posterior variance}

\hop 
The aim of the following efforts is to supplement the 
posterior mean result above with an explicit formula for 
the posterior variance of $P$, and more generally for 
the posterior variance of a $\int g\,\d P$ parameter. 
This makes construction of credibility intervals possible, 
and is used to assess full posterior consistency in the next section. 

\subsection
{\csc 6.1. Posterior variance of $P(A)$.} 
To do the posterior mean calculation, formula (5.2) sufficed. 
To calculate posterior variances requires a little list
of further formulae. Let $A_1,\ldots,A_n$ be disjoint sets, 
and let $\theta_i=P(A_i)$ with prior mean $P_0(A_i)=p_i$. 
We show later that the various means-of-products take
the following form: 
$$\eqalign{
\E\,\theta_1\theta_2\cdots\theta_n&=a_np_1\cdots p_n, \cr
\E\,\theta_1^2\theta_2\cdots\theta_n
   &=b_np_1\cdots p_n + c_np_1^2p_2\cdots p_n, \cr
\E\,\theta_1^3\theta_2\theta_3\cdots\theta_n
   &=d_np_1p_2p_3\cdots p_n+e_np_1^2p_2p_3\cdots p_n
     +f_np_1^3p_2p_3\cdots p_n, \cr
\E\,\theta_1^2\theta_2^2\cdots\theta_n
   &=g_np_1p_2\cdots p_n+h_np_1p_2(p_1+p_2)\cdots p_n
     +i_np_1^2p_2^2\cdots p_n. \cr} \eqno(6.1)$$
Here $a_n,\ldots,i_n$ are sequences of constants,
to be returned to below. 

Take a prior mean distribution $P_0$ free of atoms, and consider 
a set $A$ not meeting the data, which we again take 
to be $n$ distinct values $x_1,\ldots,x_n$. Then, 
with notation as in Section 5.2 and with $p_0=P_0(A)$,  
$$\E\,\{P(A)\midd\info_\eps\}
  ={\E\, P(A)P(A_1)\cdots P(A_n)\over \E\, P(A_1)\cdots P(A_n)}
  ={a_{n+1}p_0p_1\cdots p_n\over a_np_1\cdots p_n}
  ={a_{n+1}\over a_n}p_0, $$
while 
$$\E\,\{P(A_k)\midd\info_\eps\}
  ={\E\, P(A_k)P(A_1)\cdots P(A_n)\over \E\, P(A_1)\cdots P(A_n)}
  ={b_{n}+c_np_k\over a_n}={b_n\over a_n}+O(\eps), $$
showing that $\E\,[P\{x_k\}\midd\data]=b_n/a_n$. 
Ingredients required for second moment calculations include 
$$\E\,\{P(A)^2\midd\data\}=
  {b_{n+1}p_0p_1\cdots p_n+c_{n+1}p_0^2p_1\cdots p_n\over a_np_1\cdots p_n}
  ={b_{n+1}\over a_n}p_0+{c_{n+1}\over a_n}p_0^2, $$
$$\E\,\{P(A_k)^2\midd\info_\eps\}
  ={\E\, P(A_k)^3\prod_{i\not=k}P(A_i)\over \E\,\prod_{i=1}^nP(A_i)}
  ={d_n+e_np_k+f_np_k^2\over a_n}={d_n\over a_n}+O(\eps), $$
$$\E\,\{P(A)P(A_k)\midd\info_\eps\}
  ={b_{n+1}p_0+c_{n+1}p_0p_k\over a_n}={b_{n+1}\over a_n}p_0+O(\eps). $$
These and similar efforts entail 
$$\E\,[P\{x_k\}^2\midd\data]=d_n/a_n
  \quad {\rm and} \quad 
  \E\,[P\{x_k\}P\{x_l\}\midd\data]=g_n/a_n {\rm\ for\ }k\not=l, $$
while $\E\,[P\{x_k\}P(A)\midd\data]=(b_{n+1}/a_n)p_0$. 

Let now $A$ be any set, containing say $j$ of the data values, 
and split it into $A\cap\data=\{x_{i_1},\ldots,x_{i_j}\}$ 
and $A_0=A-\data$. Then $P(A)=P\{x_{i_1},\ldots,x_{i_{j}}\}+P(A_0)$
and, with $p_0=P_0(A)$,
$$\E\,\{P(A)\midd\data\}=(a_{n+1}/a_n)p_0 + jb_n/a_n
   =w_np_0+(1-w_n)(j/n), \eqno(6.2)$$
agreeing of course with (5.4). Next, collecting together 
the various contributions to $P(A)^2$, 
$$\E\,\{P(A)^2\midd\data\}
  =j{d_n\over a_n}+j(j-1){g_n\over a_n}
   +{b_{n+1}\over a_n}p_0+{c_{n+1}\over a_n}p_0^2 
   +2j{b_{n+1}\over a_n}p_0. \eqno(6.3)$$
We also record a formula for the cross-moment for two
disjoint sets $A$ and $B$, catching respectively 
$j$ and $k$ data points:
$$\E\,\{P(A)P(B)\midd\data\}=jk{g_n\over a_n}+j{b_{n+1}\over a_n}P_0(B)
   +k{b_{n+1}\over a_n}P_0(A)+{a_{n+2}\over a_n}P_0(A)P_0(B). $$

\subsection
{\csc 6.2. Posterior variance of a random mean.} 
We have found formulae for posterior variance of a $P(A)=\int I_A\,\d P$. 
More generally we need the posterior variance of a random mean 
$\theta=\int g\,\d P$, for which the posterior mean is given in (5.5). 
Start with a simple $g=\sum_{j=1}^my_jI_{A_j}$ with disjoint sets $A_j$,
so that $\theta=\sum_{j=1}^my_jP(A_j)$. With a little work, 
$$\eqalign{\E\,(\theta^2\midd\data)
&=\sum_{j=1}^my_j^2\Bigl\{{nd_n\over a_n}P_n(A_j)
      +{n^2g_n\over a_n}P_n(A_j)^2-{ng_n\over a_n}P_n(A_j)
      +{c_{n+1}\over a_n}P_0(A_j)^2 \cr
&\qquad\qquad  
 +{b_{n+1}\over a_n}P_0(A_j)+{2nb_{n+1}\over a_n}P_0(A_j)P_n(A_j)\Bigr\} \cr
&\qquad 
 +\sum_{j\not=k}y_jy_k\Bigl\{{n^2g_n\over a_n}P_n(A_j)P_n(A_k)
   +{nb_{n+1}\over a_n}P_n(A_j)P_0(A_k) \cr 
&\qquad\qquad +{nb_{n+1}\over a_n}P_0(A_j)P_n(A_k)
   +{a_{n+2}\over a_n}P_0(A_j)P_0(A_k)\Bigr\}, \cr}$$
which in terms of $\theta_n=\int g\,\d P_n$ and $\theta_0=\int g\,\d P_0$ 
simplifies to 
$${n^2g_n\over a_n}\theta_n^2
   +2{nb_{n+1}\over a_n}\theta_n\theta_0
   +{a_{n+2}\over a_n}\theta_0^2  
   +{b_{n+1}\over a_n}\int g^2\,\d P_0
   +{nd_n-ng_n\over a_n}\int g^2\,\d P_n. \eqno(6.4)$$
Used here is the fact that $c_n=a_{n+1}$, proved below. 

That this gives a formula $\E\,(\theta^2\midd\data)-\{\E\,(\theta\midd\data)\}^2$
for the posterior variance also for the case of any random $\int g\,\d P$,
provided only that $\int g^2\,\d P_0$ is finite, follows by 
passing to the limit via simple functions and multiple uses 
of the monotone convergence theorem. 

\subsection
{\csc 6.3. Formulae for the constants.} 
It remains to give formulae for the $a_n,\ldots,i_n$ sequences
of (6.1). We have already found that $a_n=n!\,\delta_{n-1}\cdots\delta_1$,
in terms of $\delta_j=M_{1,j}/(1-M_{0,j+1})$. For $(b_n,c_n)$,
write $\theta_1^2=_d B^2Y_1+2B\bar BY_1+\bar B^2\theta_1^2$
and $\theta_j=_d BY_j+\bar B\theta_j$ for $j=2,\ldots,n$,
as in the arguments used to prove (5.2) above. Writing out the product 
$\theta_1^2\theta_2\cdots\theta_n$ in a distributional identity
and discarding all terms involving two or more $Y_j$s 
gives an expression for its mean, which after simplification delivers 
$$\eqalign{
b_n&={1\over 1-M_{0,n+1}}\{M_{2,n-1}a_{n-1}+(n-1)M_{1,n}b_{n-1}\}, \cr
c_n&={1\over 1-M_{0,n+1}}\{2M_{1,n}a_n+(n-1)M_{1,n}c_{n-1}\}. \cr}$$
Finding $\E\,\theta_1^2$ explicitly gives start values 
$b_1=M_{2,0}/(1-M_{0,2})$ and $c_1=2M_{1,1}/(1-M_{0,2})$
for these recursive relations. Some investigations lead to 
$b_n=(n-1)!\,\delta_{n-1}\cdots\delta_1(1-w_n)$
and to $c_n=(n+1)!\,\delta_n\cdots\delta_1=a_{n+1}$. 
Working similarly with $\theta_1^3\theta_2\cdots\theta_n$ gives 
$$\eqalign{
d_n&={1\over 1-M_{0,n+2}}\{M_{3,n-1}a_{n-1}+(n-1)M_{1,n+1}d_{n-1}\}, \cr
e_n&={1\over 1-M_{0,n+2}}\{3M_{2,n}a_n+3M_{1,n+1}b_n
      +(n-1)M_{1,n+1}e_{n-1}\}, \cr
f_n&={1\over 1-M_{0,n+2}}\{3M_{1,n+1}c_n+(n-1)M_{1,n+1}f_{n-1}\}, \cr}$$
with start values 
$$d_1=M_{3,0}/(1-M_{0,3}), \quad 
  e_1=3(M_{2,1}+M_{1,2}b_1)/(1-M_{0,3}), \quad 
  f_1=3M_{1,2}c_1/(1-M_{0,3}) $$
determined from $\E\,\theta_1^3$. 
Finally, studying $\theta_1^2\theta_2^2\theta_3\cdots\theta_n$ leads to 
$$\eqalign{
g_n&={1\over 1-M_{0,n+2}}\{2M_{2,n}b_{n-1}+(n-2)M_{1,n+1}g_{n-1}\}, \cr
h_n&={1\over 1-M_{0,n+2}}\{M_{2,n}c_{n-1}+2M_{1,n+1}b_n
   +(n-2)M_{1,n+1}h_{n-1}\}, \cr
i_n&={1\over 1-M_{0,n+2}}\{4M_{1,n+1}c_n+(n-2)M_{1,n+1}i_{n-1}\} \cr}$$
for $n\ge2$, where 
$$g_2=2M_{2,2}b_1/(1-M_{0,4}), \quad  
  h_2=(M_{2,2}c_1+2M_{1,3}b_2)/(1-M_{0,4}), \quad 
  i_2=4M_{1,3}c_2/(1-M_{0,4}). $$ 
One may easily compute the $d_n,\ldots,i_n$ constants 
via these recursive schemes. 

To learn more about these sequences, observe that formula (6.3) implies 
$$1=nd_n/a_n+n(n-1)g_n/a_n+(2n+1)b_{n+1}/a_n+c_{n+1}/a_n, \eqno(6.5)$$
simply by letting $A$ be the full sample space.
Other helpful formulae for the constants involved in (6.1)
emerge as follows. Let $A_1,\ldots,A_n$ form a measurable partition 
and write $\theta_i=P(A_i)$. Then equating 
$\E\,(\sumin\theta_i)\theta_1\cdots\theta_n$ with $\E\,\theta_1\cdots\theta_n$ 
leads to $nb_n+c_n=a_n$, which with $nb_n=a_n(1-w_n)$
gives $c_n=a_{n+1}$ (again). Similarly, equating 
$\E\,(\sumin\theta_i)^2\theta_1\cdots\theta_n$ 
with $\E\,\theta_1\cdots\theta_n$ gives 
$$1={nd_n\over a_n}+{e_n\over a_n}
    +{f_n\over a_n}\sumin p_i^2+{n(n-1)g_n\over a_n}
    +{(2n-2)h_n\over a_n}+{i_n\over a_n}\sum_{i\not=j}p_ip_j. \eqno(6.6)$$
Since this is an identity valid for all $p_i$s summing to 1,
and since $\sum_{i\not=j}p_ip_j=1-\sumin p_i^2$, one must
have $f_n=i_n$ for all $n\ge 2$. Helped by this, one may show
by induction, using the recursive relations, that 
$$f_n=i_n=c_{n+1}=a_{n+2}=(n+2)!\,\delta_{n+1}\cdots\delta_1
  \quad {\rm for\ }n\ge 2. $$
Combining (6.5) with (6.6) it is also clear that 
$(2n+1)b_{n+1}/a_n+c_{n+1}/a_n=e_n/a_n+f_n/a_n+2(n-1)h_n/a_n$. 

Let us work out what happens to the iteratively defined 
sequence $d_n$. It is helpful to write 
$$d_n=(n-1)!\,y_{n-1}+(n-1)\delta_{n+1}d_{n-1} 
   \quad {\rm for\ }n\ge 1, $$
with $(n-1)!\,y_{n-1}=M_{3,n-1}(1-M_{0,n+2})^{-1}a_{n-1}$ for $n\ge1$. 
Some minutes of investigation yield 
$d_n=(n-1)!\sum_{j=0}^{n-1}y_j\delta_{j+3}\cdots\delta_{n+1}$. 
Going back to $y_n$, one sees that $y_j=\eps_j\delta_j\cdots\delta_1$,
where $\eps_j=M_{3,j}/(1-M_{0,j+2})$. Hence 
$d_n=(n-1)!\,\delta_{n+1}\cdots\delta_1\sum_{j=0}^{n-1}
   \eps_j/(\delta_{j+1}\delta_{j+2})$, 
which with $a_n=n!\,\delta_{n-1}\cdots\delta_1$ leads to 
$${nd_n\over a_n}=\delta_{n+1}\delta_n\sum_{j=1}^n
   {\eps_{j-1}\over \delta_j\delta_{j+1}}. \eqno(6.7)$$
We may similarly work out an expression for the $g_n$ sequence. Write 
$$g_n=(n-2)!\,z_{n-1}+(n-2)\delta_{n+2}g_{n-1}
  \quad {\rm for\ }n\ge2, $$
where $(n-2)!\,z_{n-1}=2M_{2,n}(1-M_{0,2+n})^{-1}b_{n-1}$;
in particular, $g_2=z_1$. This gives expressions for 
$g_3,g_4,\ldots$, and the general pattern is discovered to be 
$$g_n=(n-2)!\,(z_{n-1}+\delta_{n+1}z_{n-2}+\cdots
   +\delta_{n+1}\cdots\delta_4z_1)
  =(n-2)!\sum_{j=1}^{n-1}\delta_{n+1}\cdots\delta_{j+3}z_j. $$ 
Going back to $z_n$, which may be expressed as 
$\eta_j\delta_{j-1}\cdots\delta_1$ for 
$\eta_j=2M_{2,j+1}(1-w_j)/(1-M_{0,3+j})$, one finds 
$g_n=(n-2)!\,\delta_{n+1}\cdots\delta_1\sum_{j=1}^{n-1}
   \eta_j/(\delta_j\delta_{j+1}\delta_{j+2})$. 
In conjunction with $a_n=n!\,\delta_{n-1}\cdots\delta_1$ this implies 
$${n(n-1)g_n\over a_n}=\delta_{n+1}\delta_n\sum_{j=1}^{n-1}
   {\eta_j\over \delta_j\delta_{j+1}\delta_{j+2}}. $$ 
It will be seen in the next section that of the parts summing to 1
in (6.5) and (6.6), the $n(n-1)g_n/a_n$ is the dominant one. 

\section
\centerline{\bf 7. Consistency, and how quickly do we forget?} 

\hop 
Assume data $X_1,\ldots,X_n$ in reality follow some underlying
distribution $P_\true$. It is well known that the empirical
distribution $P_n$ converges to $P_\true$ with probability 1, 
even uniformly over all subsets, as the data volume increases. 
A question of importance is whether the Bayes estimator $\hatt P$ 
matches this feat, and, more generally, whether the posterior
distribution converges to the measure concentrated in $P_\true$. 

For parametric models it is known that Bayes inference agrees
for large samples with that based on maximum likelihood.
A more informative statement is that for Bayes and likelihood estimators 
$\hatt\theta_{B,n}$ and $\hatt\theta_{L,n}$ based on the $n$ first 
data points, it holds that 
$n^{1/2}(\hatt\theta_{B,n}-\hatt\theta_{L,n})\arr_p0$, 
even when the parametric model used to generate these likelihoods
and posteriors is incorrect, under very mild regularity assumptions;
see Hjort and Pollard (1993). It is furthermore the 
case that the posterior `forgets its prior' at a speed linear with $n$,
in the sense that aspects of the posterior traceable to the prior
has weight exactly or approximately equal to $b/(b+n)$ for a 
suitable $b$, which then can be interpreted as `prior sample size'.
The very same behaviour is observed for the Dirichlet process
prior, as shown in Ferguson (1973, 1974). We shall see that 
the situation can be quite different for other members 
of the $\GDir(H,P_0)$ class. 

\subsection
{\csc 7.1. Consistency of the posterior mean.} 
In what follows take $P_\true$ to be free of atoms on its
sample space, making all realisations $X_1,X_2,\ldots$ a.s.~distinct.
From (5.4) it is clear that $\hatt P$ also goes to $P_\true$ 
almost surely provided only that $w_n\arr0$. Under this 
key condition $\hatt P$ and the nonparametric frequentist 
estimator $P_n$ agree asymptotically. It turns out that
indeed $w_n\arr0$, but with a speed depending upon aspects 
of the distribution $H$ of the $B_j$s. 

{\smallskip\sl 
{\csc Lemma.} For any distribution $H$ for $B$, 
$w_n$ of (5.3) goes to zero with growing $n$.
\smallskip} 

{\csc Proof.} 
It suffices to show 
$$\eqalign{\E\,\bar B^{n+1}&=\int_0^1 (1-s)^{n+1}\,\d H(s)\arr 0, \cr 
(n+1)\,\E\, B\bar B^n&=\int_0^1(n+1)s(1-s)^n\,\d H(s)\arr0. \cr}$$
The first follows quickly by dominated convergence,
as does actually also the second. The point is that the 
integrand $(n+1)s(1-s)^n$ goes pointwise to zero, 
and has a maximum value bounded in $n$. Inspection shows that 
the maximum occurs for $s_0=1/(n+1)$ and that the resulting 
maximum value converges to $e^{-1}$. 
Hence there is uniform integrability and the claim follows. \square 

\smallskip 
Consider next the $\GDir(a,b,P_0)$ case, 
for which (3.5) and (5.3) yield 
$$\eqalign{w_n
&=(n+1){ab^{[n]}\over (a+b)^{[n+1]}}
  \Big/\Bigl\{1-{b^{[n+1]}\over (a+b)^{[n+1]}}\Bigr\} \cr 
&=(n+1){a\Gamma(b+n)\over \Gamma(b)}
  \Big/\Bigl\{{\Gamma(a+b+n+1)\over \Gamma(a+b)}
      -{\Gamma(b+n+1)\over \Gamma(b)}\Bigr\} \cr
&={n+1\over n+b}{a\over \Gamma(b)}
  \Big/\Bigl\{{\Gamma(a+b+n+1)\over \Gamma(b+n+1)\Gamma(a+b)}
      -{1\over \Gamma(b)}\Bigr\}. \cr} \eqno(7.1)$$
This answer generalises the well-known formula $w_n=b/(n+b)$
valid for the posterior mass outside data points for the 
Dirichlet process. Formula (7.1) gives the precise weight 
the Bayes estimator attaches to outside-of-data information, 
that is, as caused by the prior. 
The speed with which $w_n\arr0$ is different from the traditional 
$O(n^{-1})$, when $a\not=1$, as we shall see. 

Since the denominator of (5.3) goes to 1 
it suffices for large $n$ to study $u_n=(n+1)\,\E\, B\bar B^n$ 
and the speed with which this sequence tends to zero. 
For the $\GDir(a,b,P_0)$ case,
$$u_n=(n+1)a{\Gamma(a+b)\over \Gamma(b)} 
   {\Gamma(b+n)\over \Gamma(a+b+n+1)}, $$ 
and we may use the Stirling approximation, 
for example in the form of 
$$\log\Gamma(x)=(x-\half)\log x-x+\half\log(2\pi)+1/(12x)+O(1/x^2)
  \quad {\rm for\ large\ }x, $$
to assess its size. Some algebra efforts reveal  
$\log u_n=-a\log n+\log\{a\Gamma(a+b)/\Gamma(b)\}-2(a+1)+O(n^{-1})$, 
which means 
$$u_n=n^{-a}\{a\Gamma(a+b)/\Gamma(b)\}\,\exp\{-2(a+1)\}\{1+O(n^{-1})\}
  \quad {\rm when\ }n{\rm\ grows}. $$
Hence, only for the Dirichlet case $a=1$ does the posterior process
forget its origin with speed $O(n^{-1})$, which is the traditional 
speed with which memory loss sets in for Bayesian parametric statistics.
For $a>1$ the prior is forgotten more quickly and for $a<1$ more slowly
than the traditional rate.

These calculations also lead to 
$$n^{1/2}\{\hatt P(A)-P_n(A)\}\arr_p0
  \quad {\rm provided\ }a>\half. $$
Under this condition, inferential statements made by the 
Bayesian, such as credibility intervals, will agree asymptotically
with those of the frequentist using the empirical distribution. 
For smaller values of $a$, however, the speed with which the
posterior is able to forget where it came from is really too slow;
the predictive distribution is consistent, but converges slowly,
and credibility intervals will not match frequentist confidence intervals,
even for large $n$. 

\subsection
{\csc 7.2. Consistency of the posterior distribution.} 
We wish to find out whether the posterior distribution 
as such is consistent, in the sense that for any small neighbourhood
around $P_\true$, the posterior probability of such a set
converges to 1 as $n$ grows. This is a stronger statement
than merely knowing that the posterior mean is a consistent 
estimator of $P_\true$. 

{\smallskip\sl
{\csc Proposition.} 
Assume $X_1,X_2,\ldots$ are independent from some atom-free $P_\true$,
and consider $\theta=\int g\,\d P$ for an arbitrary $g$ for which 
$g^*=\int g\,\d P_\true$ is finite. Then, for almost all sample paths, 
$\theta\midd\data\arr_pg^*$. 
\smallskip} 

{\csc Proof.} 
We know that the empirical mean $\theta_n=\int g\,\d P_n$ 
goes a.s.~to $g^*$, and as above it is clear that 
$\E\,(\theta\midd\data)\arr g^*$ a.s.~in that $w_n\arr0$. 
It will suffice to show that 
$\E\,(\theta^2\midd\data)\arr (g^*)^2$ a.s.; this implies that 
the posterior variance goes to zero, and there is convergence
in probability by the usual Chebyshov inequality argument. 

To this end we work with expression (6.4), and aim to demonstrate
that $n^2g_n/a_n\arr1$ while the other terms go to zero. 
This causes $\E\,(\theta^2\midd\data)$ to go to $(g^*)^2$ 
for exactly those sample paths for which $\theta_n\arr g^*$. 
From established formulae for $a_n,b_n,c_n$ we see that 
the third and fourth terms of the right hand side of (6.5)
go to zero; this also secures that the terms $f_n/a_n$, $e_n/a_n$
and $2(n-1)h_n/a_n$ of (6.6) go to zero. It will therefore be enough 
to show that also the first term there goes to zero. 
For this we use formula (6.7). 
Note that $\delta_j\ge\delta_{j+1}$, and one finds $\eps_j\le\delta_j$.
A constant $K$ can be found such that $\delta_{j-1}/\delta_j\le K$
for all $j$. This implies 
$${nd_n\over a_n}\le\delta_{n+1}\delta_n\sum_{j=1}^n
   {\delta_{j-1}\over \delta_j\delta_{j+1}}
   \le K\delta_n\sum_{j=1}^n{\delta_{n+1}\over \delta_{j+1}}
   \le Kn\delta_n, $$ 
which goes to zero since $w_n$ does. \square 

\smallskip
Inspection of the details in these calculations show that 
the speed with which the variance goes to zero is $O(w_n)$.
As we have seen, this corresponds to the traditional
$O(1/n)$ variance rate for the Dirichlet process, whereas
the speed may be both slower and faster for the more general
prior process. 
\eject 

\section
\centerline{\bf 8. Bayesian inference and the posterior process} 

\hop
Let $P\sim{\GDir}(H,P_0)$ and assume data $x_1,\ldots,x_n$ 
have been observed. This section looks into aspects of the 
posterior process, which turns out to be quite complicated.
Only in the Dirichlet case, where $H$ is the $\Beta(1,b)$,
does the posterior seem to have an easy structure. 
Bayesian inference can nevertheless be carried out via stochastic simulation. 

\subsection 
{\csc 8.1. One data point.} 
We may take the view that $P$ of (1.2) is described in terms of 
$(B,\xi)$, where $B$ is the sequence of $B_j$s from $H$, leading
in their turn to probability weights 
$\gamma_j=\bar B_1\ldots \bar B_{j-1}B_j$, and where $\xi$ 
is the sequence of $\xi_j$s from $P_0$. Let in addition $J$ be 
a random variable in $\{1,2,3,\ldots\}$ which conditionally 
on $(B,\xi)$ has distribution given by these $\gamma_j$s,
and define $X=\xi_J$. Then $X$ given $P$ has distribution $P$. 
The task is to pass from this simultaneous representation of 
$(P,X)$ to the conditional process $P$ given $X=x$. 

When $X=\xi_J=x$ and $J=j$, one has $\xi_j=x$, without further
knowledge about the other $\xi_k$s. Furthermore, the fact that 
this happened with probability $\gamma_j$ upgrades the 
information about the distributions $B_1,\ldots,B_j$, but 
does not affect the prior information about $B_k$ for $k>j$. 
Using arguments partly parallelling those in in Sethuraman (1994, Section 4), 
one finds that 
$$P\midd\{X=x,J=j\}\sim P_{x,j}
      =\sum_{k=1}^\infty\gamma_k'\delta(\xi_k'), \eqno(8.1)$$
where on the right hand side the $\{\gamma_k'\}$ sequence 
is formed from a $\{B_k'\}$ sequence independent of the 
$\xi_k'$, which are independently drawn from $P_0$ with the exception of 
$\xi_j'$, which is equal to the fixed $x$. Now $B_k'\sim H_k'$
for $k=1,2,\ldots$, where these $H_k'$s are not equal anymore;
$\d H_k'(s)\propto(1-s)\,\d H(s)$ for $k\le j-1$,
$\d H_k'(s)\propto s\,\d H(s)$ for $k=j$, while $H_k'=H$ for $k\ge j+1$. 
Thus there is a mixture representation of the posterior as 
$$P\midd\{X=x\}\sim\sum_{j=1}^\infty q(j\midd x)P_{x,j}, 
  \quad {\rm or} \quad 
  \Pr\{P\in C\midd x\}=\sum_{j=1}^\infty q(j\midd x)\Pr\{P_{x,j}\in C\}, $$
valid for measurable subsets $C$ of the space of probability measures
on the sample space (the Borel subsets under the topology of set-wise 
convergence). It remains only to identify 
$$q(j\midd x)=\Pr\{J=j\midd X=x\}=\E\,\gamma_j=M_{0,1}^{j-1}M_{1,0}
  \quad {\rm for\ }j=1,2,\ldots. \eqno(8.2)$$
This is since the information $X=x$ from a single data point
does not change the marginal distribution $J$ has from the 
$(P,J)$ model. Notice that in (8.1) there is dependence on $x$
in the $\xi_k'$, without overburdening the notation to indicate this. 

\subsection
{\csc 8.2. The posterior in the general case.} 
Conditionally on $(B,\xi)$, the two sequences determining $P$,
let $J_1,\ldots,J_n$ be independent integer variables with
distribution given by the $\gamma_j$s, and define 
$X_1=\xi_{J_1},\ldots,X_n=\xi_{J_n}$. Then, given $P$, 
these really form an independent $n$-sample from $P$.
This provides a simultaneous representation of $(P,X_1,\ldots,X_n)$. 

Suppose for representational simplicity that the data points 
$x_1,\ldots,x_n$ are distinct. One may generalise the first result
above to 
$$P\midd\{X_1=x_1,\ldots,X_n=x_n,J_1=j_1,\ldots,J_n=j_n\}
   \sim P_{\data,j_1,\ldots,j_n}=\sum_{k=1}^\infty\gamma_k'\delta(\xi_k'), $$
where the $\{\gamma_k'\}$ is formed from a sequence of 
independent variables $\{B_k'\}$ and independently of the $\{\xi_k'\}$;
these are such that $\xi_{j_1}',\ldots,\xi_{j_n}'$ 
are fixed at values $x_1,\ldots,x_n$, respectively, while 
the remaining $\xi_k'$s are independent from $P_0$. The upgraded 
distributions $H_k'$ for $B_k'$ are given by 
$$\d H_k'(s)={\rm const.}\,(1-s)^{Y(k)-\Delta N(k)}s^{\Delta N(k)}\,\d H(s), $$
in which $Y(k)=\sum_{i=1}^n I\{j_i\ge k\}$ and 
$\Delta N(k)=\sum_{i=1}^n I\{j_i=k\}$. Hence 
$$P\midd\data\sim\sum_{j_1,\ldots,j_n{\rm\ distinct}}
  q(j_1,\ldots,j_n\midd\data)P_{\data,j_1,\ldots,j_n}. $$
It remains to give the posterior distribution of indexes.
Say that $G$ has a geometric distribution with parameter $M$ 
if $\Pr\{G=g\}=(1-M)M^g$ for $g=1,2,\ldots$. 

{\smallskip\sl
{\csc Proposition.} 
Let there be $n$ distinct data points, 
and order the random indexes $J_1,\ldots,J_n$ as 
$J_{(1)}<\cdots<J_{(n)}$. Then 
$$(J_{(1)},\ldots,J_{(n)})=_d(G_1,G_1+G_2,\ldots,G_1+\cdots+G_n), $$
where $G_1,\ldots,G_n$ are independent and geometric with 
parameters $M_{0,n},\ldots,M_{0,1}$, respectively. 
\smallskip} 

{\csc Proof.}
Knowledge of data values $\xi_{j_i}=x_i$ does not change the distribution 
of the labels as long as these are distinct. For the ordered labels 
one therefore finds the distribution 
$$\eqalign{\bar q(j_1,\ldots,j_n)
&=n!\,\E\,\gamma_{j_1}\cdots\gamma_{j_n}/\Pr(D_n) \cr 
&={n!\over \Pr(D_n)}\E\,\prod_{k=1}^\infty 
   \bar B_k^{Y(k)-\Delta N(k)}B_k^{\Delta N(k)}
 ={n!\over \Pr(D_n)}\prod_{k=1}^\infty M_{\Delta N(k),Y(k)-\Delta N(k)} \cr}$$
for $j_1<\cdots<j_n$, where $D_n$ is the event that data points are 
distinct. The product may be expressed as 
$$M_{0,n}^{j_1-1}M_{1,n-1}M_{0,n-1}^{j_2-j_1-1}M_{1,n-2}
   \cdots M_{0,2}^{j_{n-1}-j_{n-2}-1}M_{1,1}M_{0,1}^{j_n-j_{n-1}-1}M_{1,0}, $$ 
while it is shown in Section 10.2 that $\Pr(D_n)$ is equal to the
$a_n$ of formula (5.2). Combining these facts one is left with 
$$(1-M_{0,n})M_{0,n}^{j_1-1}(1-M_{0,n-1})M_{0,n-1}^{j_2-j_2-1}
  \cdots(1-M_{0,2})M_{0,2}^{j_{n-1}-j_{n-2}-1}
        (1-M_{0,1})M_{0,1}^{j_n-j_{n-1}-1}, $$
which is seen to be equivalent to the claim. \square 
\eject 

\smallskip
{\csc Remark.} 
The description above is valid for the general $\GDir(H,P_0)$ case,
and can even be generalised further to the case of different 
distributions $H_1,H_2,\ldots$ for $B_1,B_2,\ldots$ in the prior.
Note that for the particular $\GDir(a,b,P_0)$ family, 
in which the Dirichlet is the $a=1$ case, at least the $H$ to $H_k'$
updating is easy, in that 
$H_k'\sim\Beta(a+\Delta N(k),b+Y(k)-\Delta N(k))$. 
For $k$ larger than the largest $j_i$ the $H_k'$ is the same 
as the original $H$. 

For the Dirichlet case the posterior can of course be described
in a much simpler way than the scheme above. One may deduce 
from (8.1) and (8.2) that $P\midd x$ is simply another Dirichlet 
with total measure $bP_0+\delta(x)$, via various identities 
for Beta distributions; see Sethuraman (1994, Section 4). 

\section
\centerline{\bf 9. Distribution of random means}  

\hop
Recently there has been much attention given to studying
aspects of the distributions of random Dirichlet means;
see Diaconis and Kemperman (1996), Regazzini, Guglielmo and di Nunno (2000)
and Hjort and Ongaro (2000) for discussion and references. 
Here we look at the more general version of this problem,
where $P$ is a generalised Dirichlet process. 

\subsection
{\csc 9.1. General transform identities.} 
That equation (3.2) characterises the distribution of $\theta$ 
uniquely can be seen as in a parallel situation in Hjort and Ongaro (2000);
see also Lemma 3.3 in Sethuraman (1994). 
Exhibiting this distribution is however a difficult task and 
can rarely be done in closed form. The list of explicit solutions
to this problem for the Dirichlet case is very short, so a fortiori
one cannot expect explicit answers for the more general 
$\GDir(H,P_0)$ case. We point out, however, that equation (3.2)
translates into an identity for characteristic or 
moment generating functions and which can be worked with 
to extract information about the $\theta$ distribution. 
Let $L(u)=\E\,\exp(iu\theta)$ and $L_0(u)=\E_0\exp(iu Y)$. 
Via conditioning on $(B,Y)$ and then integrating over $Y$ 
one finds from (3.2) that 
$$L(u)=\int_0^1 L_0(us)L(u(1-s))\,\d H(s). \eqno(9.1)$$
In principle $L$ is determined from knowledge of $L_0$.
Similarly a convolution-type identity can be put up for the density $f$ 
of $\theta$ in terms of the density $f_0$ for $Y$ under $P_0$. 

An exception admitting a straight answer is when $Y$ is Cauchy.
One then sees that the Cauchy distribution for $\theta$ 
fits the stochastic equation (3.2), and is hence the answer; 
$\theta$ is Cauchy when $Y$ is. This is valid for any distribution 
$H$ for the $B_j$s, as can also be seen via (9.1), 
and therefore generalises a classic result for the Dirichlet process. 

\subsection
{\csc 9.2. Results for normal and stable laws.} 
Another situation of interest where some progress can be made
is the case of a normal base measure. Let $W=\sum_{j=1}^\infty\gamma_j^2$ 
in (1.2); this is a well-defined variable on $(0,1)$ with
a distribution determined via its stochastic equation 
$$W=_d B^2+(1-B)^2W, 
  \quad {\rm where\ }B\sim H{\rm\ in\ }(0,1). \eqno(9.2)$$
This follows from (1.2) in the same way as (3.1) was derived. 
If now $P_0$ is standard normal, $\theta=\sum_{j=1}^\infty\gamma_j Y_j$
is for given weights a normal $(0,W)$. This shows that $\theta$
is a scale-mixture of normals, with density of the form 
$\int_0^1\sigma^{-1}\phi(\sigma^{-1}t)p(\sigma)\,\d\sigma$, 
the $p$ density in question being the density of $W^{1/2}$. 
This density cannot be written down in closed form,
but may be arbitrarily well approximated via its moment sequence,
which may be found in a simple recursive manner; 
see Hjort and Ongaro (2000) for illustrations for the special 
Dirichlet process case. 

These arguments also work for general stable laws.
For $\alpha\in(0,2]$ and $c$ positive, say that $Y$ is stable $(\alpha,c)$ 
if its characteristic function is $\E\,\exp(iuY)=\exp(-c^\alpha|u|^\alpha)$; 
notice that $Y/c$ then is stable $(\alpha,1)$. Now take $P\sim{\GDir}(H,P_0)$
where $P_0$ is stable $(\alpha,1)$, and consider $\theta=\int x\,\d P(x)$. 
This random mean can be expressed as $\sum_{j=1}^\infty\gamma_jY_j$
where $Y_j\sim P_0$. Let $W=(\sum_{j=1}^\infty\gamma_j^\alpha)^{1/\alpha}$. 
Then $\theta$ given $\{\gamma_j\}$ is a stable $(\alpha,W)$.
It follows that $\theta$ is a scale mixture of such stable laws;
its density is $\int_0^1w^{-1}g_\alpha(w^{-1}t)\,p_\alpha(w)\,\d w$,
where $p_\alpha$ is the density of $W$ and $p_\alpha$ the 
density of a stable $(\alpha,1)$ variable. 

\section
\centerline{\bf 10. Concluding remarks} 

\hop
In these final remarks a couple of further uses of the 
generalised Dirichlet process are identified, and possibilities
for further research are noted. 

\smallskip
{\csc 10.1. Bayesian robustness.} 
If a statistician uses the Dirichlet $(b_0,P_0)$ process as a prior,
or as an element in a more complicated prior, one may
supplement such analysis with that using the $\GDir(a,b,P_0)$
prior, preferably with the proviso $2b/(1+a)=b_0$,
as indicated in Section 4. Answers derived under the Dirichlet
should then be compared to those obtained with the more general
prior, say corresponding to values of $a$ inside $(\half,{3\over 2})$.
Small variation in results indicates Bayesian robustness. 

\smallskip
{\csc 10.2. Marginal distribution when data are distinct.} 
Let $P\sim{\GDir}(H,P_0)$ with consequent observations $X_1,X_2,\ldots$. 
Consider $D_n$, the event that the $n$ first observations are distinct. 
From the definition (1.2), 
$$\Pr(D_n)=n!\,\Pr\{X_1<\cdots<X_n\}
   =n!\sum_{i_1<\cdots<i_n}\E\,\gamma_{i_1}\cdots\gamma_{i_n}. $$
But from the proof of the lemma of Section 5.1 it is clear that 
$\Pr(D_n)=a_n=n!\,\delta_{n-1}\cdots\delta_1$, in the notation
of Section 5.3. It also follows that for disjoint sets $A_1,\ldots,A_n$, 
$$\Pr\{X_1\in A_1,\ldots,X_n\in A_n\midd D_n\}=P_0(A_1)\cdots P_0(A_n), $$
that is, conditional on data points being distinct, 
the observations form an i.i.d.~sequence from $P_0$.
This generalises a result for the Dirichlet process due to 
Korwar and Hollander (1973).

\smallskip
{\csc 10.3. A semiparametric prior giving density estimates.} 
Assume that $Y_i=\theta+\eps_i$ for $i=1,\ldots,n$ with 
$\eps_1,\ldots,\eps_n$ being independent from a $P$ 
centred at zero. For this signal plus noise model a sensible prior 
could be to give $\theta$ a prior $\pi(\theta)\,\d\theta$ 
and $P$ an independent $\GDir(H,P_0)$ process prior,
where $P_0$ has a density $p_0$ centred at zero.  
Then calculations similar to but more general than those of Section 5
show that $\theta$ given observations $y_1,\ldots,y_n$ 
has posterior density 
$\pi(\theta\midd\data)=c\,\pi(\theta)\prod_{i=1}^n p_0(y_i-\theta)$,
which is also the posterior computed under the simple parametric model
where $P=P_0$. It is assumed here that the $y_i$s are distinct. 
Since knowing data and $\theta$ amounts to knowing the $\eps_i$s, 
results of Sections 5 and 6 apply, giving 
$$\E\,\{P(A)\midd\data,\theta\}=w_nP_0(A)
   +(1-w_n)n^{-1}\sumin I\{y_i-\theta\in A\}. $$
But this gives
$$\hatt P(A)=\E\,\{P(A)\midd\data\}
   =w_nP_0(A)+(1-w_n)n^{-1}\sumin \Pr\{\theta\in y_i-A\midd\data\}, $$
which is found to be an integral of a smooth density estimate, 
$$\hatt p(t)=w_np_0(t)+(1-w_n)n^{-1}\sumin\pi(y_i-t\midd\data). $$
This is a mixture of the prior guess density and a kernel type
density estimator, with bandwidth approximately proportional to $n^{-1/2}$.
The construction here can be generalised to include scale parameters
and covariates. 

\smallskip
{\csc 10.4. Prior process with different $B_i$ distributions.} 
As the complicated posterior indicates, it may be useful
to allow different distributions $H_1,H_2,\ldots$ 
for the $B_1,B_2,\ldots$ in (1.2). A condition guaranteeing 
a.s.~convergence of $\bar B_1\cdots\bar B_n$ to zero is needed. 
Tsilevich (1997) has actually worked with a particular construction 
of this type, but in a different probabilistic framework, 
and she does not discuss applications or implications for Bayesian statistics. 
For the general prior process indexed by $H_1,H_2,\ldots$ 
the posterior of $P$ given a set of data becomes of the same type,
with updated $H_1',H_2',\ldots$, following the lines of Section 8.
Accordingly, at least in a technical sense of the term, 
we have constructed a large conjugate class of nonparametric priors. 

\smallskip
{\csc 10.5. Ties in data.} 
Formulae for posterior mean and variance were derived above 
for the case of data points $x_1,\ldots,x_n$ being distinct,
as they would be if stemming from an underlying atom-free distribution.
When the $X_i$s really come from a $P$ chosen by the generalised
Dirichlet process there will be multiple ties with positive probability,
however. A more complete description should therefore include 
generalised versions of say (5.4) and (6.4) for multiplicities 
among the data points. This is possible but requires cumbersome 
extensions of arguments and recursive schemes developed in Section 6.3. 
To illustrate, and to compare issues of data weighting with 
the distinct case and with the Dirichlet case, we indicate here 
results for the case of $x_1=x_2$ distinct from $n-2$ distinct values 
$x_3,\ldots,x_n$. Let $\info_\eps$ indicate the information 
$X_1,X_2\in S(x_1,\eps)$ and $X_i\in S(x_i,\eps)$ for $i=3,\ldots,n$, 
and write $\theta_i=P(S(x_i,\eps))$. With arguments and notation 
as in Section 6.1 one first finds that 
$$\E\,\{P(A)\midd\info_\eps\}
   ={\E\, P(A)\theta_1^2\theta_3\cdots\theta_n
    \over \E\,\theta_1^2\theta_3\cdots\theta_n}
   ={b_nP_0(A)+c_nP_0(A)p_1\over b_{n-1}+c_{n-1}p_1}
   ={b_n\over b_{n-1}}P_0(A)+O(\eps) $$
for sets $A$ not meeting the data, which means that 
$\hatt P=\E\,\{P(\cdot)\midd\data\}$ is the same as $(b_n/b_{n-1})P_0$
outside the data set. Furthermore, 
$$\E\,\{P(A_1)\midd\info_\eps\}
   ={\E\,\theta_1^3\theta_3\cdots\theta_n
    \over \E\,\theta_1^2\theta_3\cdots\theta_n}
   ={d_{n-1}+e_{n-1}p_1+f_{n-1}p_1^2\over b_{n-1}+c_{n-1}p_1}
   ={d_{n-1}\over b_{n-1}}+O(\eps), $$
$$\E\,\{P(A_3)\midd\info_\eps\}
   ={\E\,\theta_1^2\theta_3^2\cdots\theta_n
    \over \E\,\theta_1^2\theta_3\cdots\theta_n}
   ={g_{n-1}+h_{n-1}(p_1+p_3)+i_{n-1}p_1p_3
    \over b_{n-1}+c_{n-1}p_1}
   ={g_{n-1}\over b_{n-1}}+O(\eps). $$
Accordingly, $\E\,[P\{x_1\}\midd\data]=d_{n-1}/b_{n-1}$
while $\E\,[P\{x_i\}\midd\data]=g_{n-1}/b_{n-1}$ for the $n-2$
other data points. 

For the $\GDir(a,b,P_0)$ process one learns that for $a<1$, 
there is slightly less weight $b_n/b_{n-1}$ to the outside-the-data
set with the $x_1=x_2$ tie than without such a tie;
the situation is reversed for $a>1$. The expected weight $d_{n-1}/b_{n-1}$
given to the double data point $x_1$ can similarly be compared with 
$2b_n/a_n$, the expected weight given to $\{x_1,x_2\}$ when these
are distinct. Here sometimes the first is bigger than the second
and sometimes the other way around, for a given $a\not=1$. 
Quite generally, these probability weights given to the outside-the-data
set and to the individual data points are independent of ties
if and only if the process is a Dirichlet, that is, the 
$H$ distribution is a $\Beta(1,b)$ for some $b$. 

\smallskip
{\csc 10.6. Other nonparametric priors.} 
One sees from the results of Section 8 that the posterior distribution
of an arbitrary random mean $\theta=\int g\,\d P$ has the structure 
$$\theta\midd\data\sim D_0T+\sumin D_ig(x_i), \eqno(10.1)$$
where $D_0,D_1,\ldots,D_n$ are random weights summing to 1 
and $T$ is a variable with mean $\theta_0=\int g\,\d P_0$. The distribution 
of $(D_1,\ldots,D_n)$ is symmetric when the data values are distinct,
securing equal weight to each data point. The distribution is 
more complicated with ties in the data, as indicated above. 
It is interesting that several different unrelated nonparametric priors 
lead to the structure (10.1), among them two constructions of 
Hjort and Ongaro (2001), and, of course, the Dirichlet. For each
such prior the predictive distribution takes the form 
$w_nF_0+(1-w_n)\tilda F_n$, say, where $w_n=\E\, D_0$ and 
$\tilda F_n$ is a `modified empirical distribution function',
being equal to the empirical $F_n$ when data are distinct and otherwise
awarding somewhat modified weights to data points 
with different multiplicities. A characterisation theorem of Lo (1991) 
implies that only when the prior is a Dirichlet do these weights 
become proportional to the multiplicities, that is, only then 
is $\tilda F_n$ equal to $F_n$ for all data configurations. 

\section 
\centerline{\bf References}

\def\ref#1{{\noindent\hangafter=1\hangindent=20pt
  #1\smallskip}}          
\parindent0pt
\baselineskip11pt
\parskip3pt 
\medskip 

\ref{%
Diaconis, P.~and Kemperman, J. (1996). 
Some new tools for Dirichlet priors.
In {\sl Bayesian Statistics 5}, 
J.M.~Bernardo, J.O.~Berger, A.P.~Dawid and A.F.M.~Smith (eds.), 97--106. 
Oxford University Press, Oxford.} 

\ref{%
Feigin, P.B.~and Tweedie, R.L. (1989).
Linear functionals and Markov chains associated with Dirichlet processes. 
{\sl Mathematics Proceedings of the Cambridge Philosophical Society}
{\bf 105}, 579--585.}

\ref{%
Ferguson, T.S. (1973).
A Bayesian analysis of some nonparametric problems.
{\sl Annals of Statistics} {\bf 1}, 209--230.}

\ref{%
Ferguson, T.S. (1974).
Prior distributions on spaces of probability measures.
{\sl Annals of Statistics} {\bf 2}, 615--629.}

\ref{%
Guglielmi, A., Holmes, C.C.~and Walker, S.G. (2001).
Perfect simulation involving functionals of a Dirichlet process.
To appear.} 

\ref{%
Guglielmi, A.~and Tweedie, R.L. (2000).
MCMC estimation of the law of the mean of a Dirichlet process.
Technical report TR 00.15, CNR--IAMI, Milano.} 

\ref{%
Hjort, N.L. (2001). 
Topics in nonparametric Bayesian statistics.
In {\sl Highly Structured Stochastic Systems},
to be published.} 

\ref{%
Hjort, N.L.~and Ongaro, A. (2000). 
On the distribution of random Dirichlet means. 
Statistical Research Report, University of Oslo.} 

\ref{%
Hjort, N.L.~and Ongaro, A. (2001). 
Two generalisations of the Dirichlet process.
In progress.}

\ref{%
Hjort, N.L.~and Pollard, D. (1993). 
Asymptotics for minimisers of convex processes.
Statistical Research Report, University of Oslo.} 

\ref{%
Korwar, R.M.~and Hollander, M. (1973).
Contributions to the theory of Dirichlet processes.
{\sl Annals of Probability} {\bf 1}, 705--711.} 

\ref{%
Lo, A.Y. (1991).
A characterization of the Dirichlet process. 
{\sl Statistics and Probability Letters} {\bf 12}, 185--187.} 

\ref{%
Paulsen, J.~and Hove, A. (1999).
Markov chain Monte Carlo simulation of the distribution
of some perpetuities.
{\sl Advances of Applied Probability} {\bf 31}, 112--134.}

\ref{%
Regazzini, E., Guglielmi, A.~and di Nunno, G. (2000).
Theory and numerical analysis for exact distributions
of functionals of a Dirichlet process.
Research report, Universit\`a di Pavia.} 

\ref{%
Sethuraman, J. (1994).
A constructive definition of Dirichlet priors.
{\sl Statistica Sinica} {\bf 4}, 639--650.} 

\ref{%
Sethuraman, J.~and Tiwari, R. (1982).
Convergence of Dirichlet measures and the interpretation of their
parameter. In {\sl Proceedings of the Third
Purdue Symposium on Statistical Decision Theory and Related Topics},
S.S.~Gupta and J.~Berger (eds.), 305--315. Academic Press, New York.}

\ref{%
{\cyr Tsilevich, N.V. (1997).
Raspredelenie srednego znacheniya dlya nekotoryh slu\-cha{\j3}nyh mer.
Zapiski nauchnyh seminarov POMI}, 
Petersburg Department of Mathematical Institute, 
University of Sankt-Peterburg, 268--279.} 

\ref{%
Walker, S.G., Damien, P., Laud, P.W.~and Smith, A.F.M. (1998).
Bayesian nonparametric inference for random distributions and related 
functions (with discussion). 
{\sl Journal of the Royal Statistical Society} {\bf B 61}, 485--527.} 

\bye